\documentclass[oneside,english]{amsart}
\usepackage[T1]{fontenc}
\usepackage[latin9]{inputenc}
\usepackage{amsthm}
\usepackage{amssymb}
\usepackage{graphicx}

\numberwithin{equation}{section}
\numberwithin{figure}{section}

\usepackage{babel}

\usepackage{babel}

\makeatother

\usepackage{babel}
\begin{document}

\vspace*{2cm} \normalsize \centerline{\Large \bf Singularly Perturbed Profiles}

\vspace*{1cm}

\centerline{\bf V.Bykov$^{a}$, Y.Cherkinsky, V.Gol'dshtein, N.Krapivnik$^{b,*}$\footnote{Corresponding
author. E-mail: mordeev@bgu.ac.il}, U.Maas$^{a}$ }

\vspace*{0.5cm}

\centerline{$^{a}$ Institute of Technical Thermodynamics, Karlsruhe University(TH), Kaiserst.12, 76128 Karlsruhe, Germany}

\centerline{$^{b}$ Department of Mathematics, Ben Gurion University of the Negev,
P.O.B. 653, 84105 Beer-Sheva, Israel}


\vspace*{1cm}

\noindent {\bf Abstract.}

In the current paper the so-called REaction-DIffusion Manifold (REDIM)
method of model reduction is discussed within the framework of standard
singular perturbation theory. According to the REDIM a reduced model
for the system describing a reacting flow (accounting for chemical
reaction, advection and molecular diffusion) is represented by a low-dimensional
manifold, which is embedded in the system state space and approximates
the evolution of the system solution profiles in space and in time.
This pure geometric construction is reviewed by using Singular Perturbed
System (SPS) theory as the only possibility to formalize, to justify
and to verify the suggested methodology. The REDIM is studied as a
correction by the diffusion of the slow invariant manifold defined
for a pure homogeneous system. A main result of the study is an estimation
of this correction to the slow invariant manifold. A benchmark model
of Mechaelis-Menten is extended to the system with the standard diffusion
described by the Laplacian and used as an illustration and for validation
of analytic results. 


\section{Introduction}

Recently, computational tools of numerical integration have become
a prevailing methodology in study of various phenomena in reacting
flow systems. For instance, nowadays in combustion theory not only
the development of practical devices, which utilize combustion of
hydrocarbon fuels, but also theoretical research in the field to a
great extend are based on the numerical integration of the mathematical
models with detailed kinetics \cite{WMD2004,B2008}. The problems
of pollutant mitigation \cite{DGA2009}, increase of efficiency, safety,
controllability and optimization of combustion facilities and processes
cannot be handled without numerical computations using detailed knowledge
on the kinetics of combustion \cite{TL2009,SW2014}. To date there
is a number of numerical packages where the detailed reaction mechanisms
of oxidation of various hydrocarbon fuels are implemented (see e.g. \cite{chemkin}).
There is also a number of very detailed kinetic mechanisms developed specifically
to describe the high temperature combustion processes of hydrocarbon
oxidation, where all possible reaction paths are accounted for in
order to describe combustion processes reliably. Similar situation
in other fields dealing with reacting flows, namely, the kinetic models
constantly grow in complexity (non-linearity) and in dimension, e.g.
in systems biology, atmospheric chemistry etc.

There are two principal ways to overcome overwhelming growth of mechanisms
of chemical kinetics. The first is an accurate modeling focusing on
a particular application. Such models can be very accurate, but typically
they are valid in very narrow range of the system physical parameters.
At this level standard methodologies, as for instance Quasi-Steady
State Assumption (QSSA) and Partial Equilibrium Assumption (PEA) (see
e.g. \cite{WMD2004,Echekki2011}), can be applied to obtain very compact descriptions
even analytically. However, extrapolation of results obtained by using
such models might be in question.

The second more attractive alternative is an automatic model reduction
approaches aiming at a compact formulation of a reduced model \cite{OM1998},
which is almost completely congruent to the original detailed model.
In contrast to conventional methodologies of model reduction, e.g.
QSSA and PEA where typically the number of elementary reactions (steps)
and species (dimension) is reduced, the methods of automatic model
reduction are based on the low-dimensional manifolds (e.g. ILDM, FGM,
REDIM etc. \cite{MP1992,P1987,OG2000,BM2007}) in which though the
dimension is reduced no elementary reaction/specie is neglected.
In a number of works (e.g. \cite{BM2007,BG2008,BGM2008,BM2009,BM2009ZFC,BNKM2015})
the progress in the development of this methodology was reported.
These methods fully preserve the topological and dynamical structure
of chemical reaction networks. They are based on so-called invariant
manifolds of low-dimension, which define constrains (relations) between
the variables of the detailed model and therefore can be efficiently
used to formulate a consistent (with the overall dynamics of the detailed
model) reduced model \cite{B2011,BGM2011,MB2011}.

In the suggested study we focus on one of such methods for the reduced
model formulation, namely, on the so-called REaction-DIffusion Manifold
(REDIM) approach \cite{BM2007}. In this approach the constrains,
which define a manifold such that 
\begin{itemize}
\item manifolds are slow (i.e. based on the concept of decomposition of
motions, which reduces also the stiffness of the reduced model); 
\item attractive (i.e. guaranties stability with respect to small perturbations); 
\item invariant (i.e. as itself represents the system solution of the detailed
model, which guaranties the accuracy, consistency and congruence of
the reduced model); 
\item low-dimensional (i.e. guarantees dimension reduction). 
\end{itemize}

After the manifold with needed properties has been identified,it can be used to
reduce the system by exploiting reformulation of the system of governing
equations for the independent set of system variables (parameters
of a manifold) and use them for modelling of reacting flow. As a result
we obtain a reduced system \textendash{} reduced model. The properties
above guarantee that not only topology of the original mechanism of
chemical kinetics \textendash{} chemical reaction network will be
preserved, but the information on the species concentrations/variables
will not be lost by reduction implemented in this way, because the
manifold gives us possibility to compute the dependent variables by
using independent ones and, therefore, we obtain the information about
the whole thermo-chemical state space of reacting systems.

In the next section the theory of slow invariant manifolds for a homogeneous
system is briefly outlined for completeness of the exposition. Then,
the REDIM method will be discussed as a method to construct the manifold
approximating the evolution of the detailed system solution profiles.
Thus, the SPS system with an additional Laplacian operator is taken
to represent mathematically simplest possible model of the REDIM.
Afterwords the standard SPS theory is applied to analyze the role
of the transport (the limit of the slow transport is considered) with
the first and second correction terms to the slow invariant manifold
obtained implicitly. Finally, main results are illustrated by Michaelis-Menten
model with diffusion.

\section{Slow invariant manifolds}

We start from an intuitive definition of a smooth manifold that belongs
to an open subset $U$ of the Euclidean space $R^{m}$, which defines
the reacting system state space. An $d$ - dimensional subset $M$
of $U$ is called an $d$ - dimensional manifold if each point $x\in M$
has a neighborhood $V_{x}\subset M$ that is diffeomorphic to the
$d$ - dimensional unit ball $B^{d}$ (i.e there exist a smooth map
$\varphi:V_{x}\to B^{d}$ with the smooth inverse map). The tangent
space at the point $x\in M$ is a collection of all vectors (``directions'')
tangent to $M$ at $x$. The tangent space is a linear $d$ - dimensional
subspace of $R^{m}$.

Let $F(x)$ be a smooth vector field defined on $U$. A smooth manifold
$M$ is an invariant manifold for $F(x)$ if $F(x)$ is a tangent
vector field to $M$, i.e. $F(x)\in M$. As a result any trajectory $u(t)$
of the system 
\begin{equation}
\frac{du}{dt}=F(u)\label{eq:inv}
\end{equation}
that has an initial point $u_{0}$ in $M$. If $M\in R^{m}$ and has not a boundary, then
$u(t)$ belongs to $M$ for any $t$. An alternative description that
do not use initial data is the following: A smooth manifold $M$ is
an invariant manifold for $F(x)$ if a trajectory of (\ref{eq:inv})
has a common point $u(t_{0})$ with $M$ then it belongs to $M$ for
$t$ close to $t_{0}.$

Consider the system of ordinary differential equations in which derivatives
of a number of variables are multiplied by a small parameter i.e.
Singularly Perturbed System (SPS). The conventional form of the SPS
is given by 
\begin{equation}
\frac{dx}{dt}=F_{s}(x,y,\epsilon)\label{eq:sps1}
\end{equation}

\begin{equation}
\epsilon\frac{dy}{dt}=F_{f}(x,y,\epsilon)\label{eq:sps2}
\end{equation}

Here $x\in\mathbb{R}^{n_s}$, $y\in\mathbb{R}^{m_f}$ are vectors in Euclidean
spaces, $t\in(t_{0},\infty)$ is a variable with the meaning of time
and $0<\epsilon<\epsilon_{0}\ll1$ is a small positive parameter. In
all cases of interest here, the functions $F_{s}:\mathbb{R}^{n_s}\times\mathbb{R}^{m_f}\longrightarrow\mathbb{R}^{n}$,
$F_{f}:\mathbb{R}^{n_s}\times\mathbb{R}^{m_f}\longrightarrow\mathbb{R}^{m}$
are infinitely differentiable for all $x\in\mathbb{R}^{n_s}$, $y\in\mathbb{R}^{m_f}$
, $0<\epsilon<\epsilon_{0}$ (at least, in the relevant domain). The
values $|f_{i}(x,y,\epsilon)|$ and $|g_{j}(x,y,\epsilon)|$ ($i=1,2,\ldots n$;
$j=1,2,\ldots m$) are assumed to be comparable with unity as $\epsilon\rightarrow0$.

The vector field corresponding to the system (\ref{eq:sps1}-\ref{eq:sps2})
is 
\[
F_{\epsilon}(x,y)=(F_{s}(x,y,\epsilon),\frac{1}{\epsilon}F_{f}(x,y,\epsilon))
\]

A usual approach for a qualitative study of SPS systems is to consider
first the algebraic differential system (so-called degenerate system)
\[
\frac{dx}{dt}=F_{s}(x,y,\epsilon)
\]

\[
0=F_{f}(x,y,0)
\]
and then to draw conclusions about the qualitative behavior of the
full system for sufficiently small $\epsilon$.

In order to recall a basic notions of the theory of singularly perturbed
systems, a system of equations 
\[
\frac{dx}{dt}=F_{s}(x,y,\epsilon)
\]
is called the slow subsystem, $x$ is called the slow variable and
the system of equations 
\[
\epsilon\frac{dy}{dt}=F_{f}(x,y,\epsilon)
\]
is called the fast subsystem, $y$ is called the fast variable.

Equation 
\[
0=F_{f}(x,y,0)
\]
determines the slow manifold $S$ of the system (\ref{eq:sps1}-\ref{eq:sps2})
\cite{BGG2006,BG2008}. Points of the slow manifold are sub-divided
into two types, standard points and turning points. A point $(x,y)$
is a standard point of the slow manifold if in some neighborhood of
this point the slow manifold $S$ can be represented as a graph of
the function $y=h_{0}(x)$ such that $F_{f}(x,h_{0}(x),0)=0$. Practically
it means that the condition of the Implicit Function Theorem $det(D_{y}F_{f}(x,h_{0}(x),0))\neq0$
holds and the slow manifold has the dimension of the slow variable
$x$. Points where this condition does not hold are called turning
points of the slow manifold.

We are looking for a family $M_{\epsilon}$ of $n$-dimensional smooth
invariant manifolds of the vector field
$F_{\epsilon}(x,y)=(F_{s}(x,y,\epsilon),\frac{1}{\epsilon}F_{f}(x,y,\epsilon))$
that can be represented as a graph of the function $y=h(x,\epsilon)$
and satisfy an additional property $\lim_{\epsilon\to0}h(x,\epsilon)=h_{0}(x)$.
Here $y=h_{0}(x)$ is an analytic representation of a sheet of $S$.
For any fixed $\epsilon$ the invariant manifold $M_{\epsilon}$ is
called the slow invariant manifold or the manifold of slow motions.
[[For $\epsilon=0$ the invariant manifold $M_{0}$ coincides with slow manifold $S$, in other words,
slow manifold $S$ is zero ($\epsilon=0$) approximation of the slow invariant
 manifolds $M_{\epsilon}$. ]]
Under reasonable analytic assumptions was proved existence of $M_{\epsilon}$
for any comparatively small $\epsilon$.

The motion along the slow manifold is described by the equation: 
\[
\frac{dx}{dt}=F_{s}(x,h_{0}(x),0),
\]
and the motion along the slow invariant manifold is described by the
equation: 
\[
\frac{dx}{dt}=F_{s}(x,h(x,\epsilon)),\epsilon).
\]

In general situations, the determination of the exact form and location
of the slow invariant manifold is impossible. Therefore, methods of
approximation are necessary. One of them finds the slow invariant
manifold as a power series with respect to the small parameter $\epsilon$
\[
h(x,\epsilon)=h_{0}(x)+\sum_{i=1}^{\infty}\epsilon^{i}h_{i}(x).
\]

It is clear that the slow manifold $y=h_{0}(x)$ is an $O(\epsilon)$
 approximation of the slow invariant manifolds. Thus, the general
scheme of application of this technique for a singularly perturbed
system can be subdivided into the analysis of the fast and slow motions.
The analysis can be considerably simplified by this decomposition
and by reducing the dimension of the system to the dimension of the
slow variable $x$ and the dimension of the fast variable $y$.

Note that the formal substitution of the function $h(x,\epsilon)$
instead $y$ into (\ref{eq:sps2}) gives the first order PDE, the
so-called invariance equation 
\[
\epsilon\frac{\partial h}{\partial x}(x,\epsilon)F_{s}(x,h(x,\epsilon)),\epsilon)=F_{f}(x,h(x,\epsilon)),\epsilon).
\]
for $h(x(t),\epsilon)$, since $\epsilon\frac{dy}{dt}=\epsilon\frac{\partial h}{\partial x}\frac{dx}{dt}$.

However, it is generally not possible to find the explicit solution
$y=h_{0}(x)$ exactly from the equation $0=F_{f}(x,y,0)$. In this
case the slow invariant manifold may be obtained in an implicit form
\[
G(x,y,\epsilon)=0
\]

(that satisfies the following property: $G(x,y,0)=F_f(x,y,0)$)
or, at least, the function $G(x,y,\epsilon)$ can be approximated
in the implicit form.

The implicit form of the invariance equation is 
\[
\frac{\partial G}{\partial y}F_{f}(x,y,\epsilon)+\epsilon\frac{\partial G}{\partial x}F_{s}(x,y,\epsilon)=0.
\]

The verification is based on the Implicit Function Theorem.

To obtain the first order approximation in the implicit form , it
is necessary to differentiate the equation $F_{f}(x(t),y(t),0)=0$
as a function of $t$ and by virtue of the system (\ref{eq:sps1}-\ref{eq:sps2}),
the result is 
\[
\frac{\partial F_{f}(x,y,0)}{\partial y}F_{f}(x,y,\epsilon)+\epsilon\frac{\partial F_{f}(x,y,0)}
{\partial x}F_{s}(x,y,\epsilon)=0.
\]

If $\det\frac{\partial F_{f}}{\partial y}\neq0$ this equation can
be written in the more convenient form 
\[
F_{f}(x,y,\epsilon)+\epsilon\left[\frac{\partial F_{f}}{\partial y}\right]^{-1}\frac{\partial F_{f}}
{\partial x}F_{s}(x,y,\epsilon)=0.
\]

We recover the zero order approximation, on setting $\epsilon=0$.
To obtain the second order approximation, it is necessary to differentiate
the equation $F_{f}(x(t),y(t),0)=0$ twice in the respect of
$t$. Suppose $\det\frac{\partial F_{f}}{\partial y}\neq0$. Let
us use the short notation $N:=\left[\frac{\partial F_{f}}{\partial y}\right]^{-1}\frac{\partial F_{f}}{\partial x}F_{s}(x,y,\epsilon)$.
After elementary calculations we have 
\[
F_{f}(x,y,\epsilon)+\epsilon\left[N+\left[\frac{\partial F_{f}}{\partial y}\right]^{-1}\frac{\partial N}{\partial y}F_{f}(x,y,\epsilon)\right]+\epsilon^{2}\left[\frac{\partial F_{f}}{\partial y}\right]^{-1}\frac{\partial N}{\partial x}F_{s}(x,y,\epsilon)=0.
\]

This is the second approximation of the slow invariant manifold in
the implicit form.

\subsection{Singularly perturbed systems (SPS) and vector fields (SPVF)}

A singularly perturbed system (fast-slow system ) of ordinary differential
equations is considered as the main mathematical construction to handle
the multiple time scales \cite{F1979}, \cite {GS1992}. Below the standard
SPS is given in the fast time (renormalization of time $\tau=\varepsilon t$
is implemented) formulation to underline the main idea behind the
framework of singularly perturbed vector fields (see \cite{BGG2006}
for more details): 
\begin{equation}
\begin{array}{c}
\frac{dx}{d\tau}=F_{f}\left(x,\,y\right)\\
\frac{dy}{d\tau}=\varepsilon\,F_{s}\left(x,\,y\right)
\end{array}\label{eq:SPS}
\end{equation}
Here $\varepsilon\ll1$ is a small parameter, $x=\left(x_{1},...,\,x_{m_{f}}\right)$
is the fast vector variable, $y=\left(y_{1},...,\,y_{m_{s}}\right)$
is the slow variable, $n=m_{s}+m_{f}$. The time variable is usual
for models without renormalization. When the system is considered
in the general form 
and difference in characteristic time scales is assumed, then such
system can be treated as the system of ordinary differential equations
with a small parameter $\delta$ where a priori a division on fast
and slow sub-sytems is not known, but the vector field $F(u,\delta)$
has a ``hidden'' fast-slow structure. It means that the vector variable
$u$ has no explicit fast and slow components. Thus, the main idea
of SPVF is that after a suitable coordinate transformation $u\mapsto\left(x,\,y\right)$
the original system (\ref{eq:inv}) can be transferred to the standard
singularly perturbed system (\ref{eq:SPS}) with a new small parameter,
$\varepsilon$ that depends on $\delta$ (see the appendix for more
details). After the transformation is , one can apply all powerful
technique of singular perturbations approximations fast ($y=const$)
and slow manifolds ($F_{f}\left(x,\,y\right)=0$) etc. Then, by using
 transformation one can trace the information (decomposition, fast
and slow manifolds and their properties etc.) the original variables
$u$.

More details on the concept and the theory of singular perturbed vector
fields (SPVF) can be found in \cite{BGG2006, BGM2008} . For the implementation,
however, two main problems have to be solved:

1) How to check that a vector field under consideration is a singularly
perturbed one?

2) How to find a corresponding coordinate transformation that transforms
the original system to the standard SPS system?

At the moment,  to both questions in the general case , but in the
case of a linear transformation of coordinates corresponding algorithm
called as a global quasi-linearization (GQL) has been developed (see
e.g. \cite{BG2013, BGM2008, BM2009_a}).

\subsection{GQL and system decomposition}

In the case when fast manifolds and the system decomposition have
linear structure they can be identified by a gap between the eigenvalues
of an appropriate global linear approximation of the Right Hand Side
(RHS) - vector function of (\ref{eq:inv}) (see \cite{BM2009} for
detailed discussion)

\[
Tu\approx F\left(u\right).
\]

Note that we did not use a hidden small parameter $\delta$ in $F\left(u\right)$,
because its existence is not known a priori and has to be validated
in a course of application of the GQL. Now, if $T$ has two groups
of eigenvalues: the so-called small eigenvalues $\lambda\left(\Lambda_{s}\right)$
and large eigenvalues $\lambda\left(\Lambda_{f}\right)$  have sufficiently
different order of magnitude, then the vector field $F(u)$ is
regarded as linearly decomposed asymptotic singularly perturbed vector
field \cite{BGG2006}. Accordingly, fast and slow invariant subspaces
given by columns of the matrices $Z_{f},\:Z_{s}$ corresponding 
\cite{MP1992a} define the slow and  variables . Namely, 

\begin{equation}
T\equiv\begin{pmatrix}Z_{f} & Z_{s}\end{pmatrix}\cdot\begin{pmatrix}\Lambda_{f} & 0\\
0 & \Lambda_{s}
\end{pmatrix}\cdot\begin{pmatrix}\tilde{Z_{f}}\\
\tilde{Z_{s}}
\end{pmatrix},\label{eq:11_T_GQL}
\end{equation}

now, if we denote 

\[
\tilde{Z}=Z^{-1}=\begin{pmatrix}Z_{f} & Z_{s}\end{pmatrix}^{-1}=\begin{pmatrix}\tilde{\left(Z_{f}\right)}_{m_{s}\times n}\\
\left(\tilde{Z_{s}}\right)_{m_{f}\times n}
\end{pmatrix},
\]

then, new coordinates suitable for an explicit decomposition (and
coordinates transformation) are given by $(U,\,V)$:

\begin{equation}
\begin{array}{c}
U:=\tilde{Z}_{f}\,u\\
V:=\tilde{Z}_{s}\,u
\end{array}.\label{eq:13_VARIABLES-2}
\end{equation}

The decomposed form and corresponding fast and slow subsystems becomes

\begin{equation}
\left\{ \begin{array}{c}
\frac{dU}{dt}=\tilde{Z}_{f}\cdot F\left(\left(Z_{f}\:Z_{s}\right)\left(\begin{array}{c}
U\\
V
\end{array}\right)\right)\\
\frac{dV}{dt}=\tilde{Z}_{s}\cdot F\left(\left(Z_{f}\:Z_{s}\right)\left(\begin{array}{c}
U\\
V
\end{array}\right)\right)
\end{array}\right..\label{eq:14_ODE_DECOMPOSED}
\end{equation}

The small system parameter controlling the characteristic time scales
in (\ref{eq:14_ODE_DECOMPOSED}) can be estimated by the gap between
the smallest eigenvalue of the large group and the largest eigenvalue
of the small group of eigenvalues \cite{BM2009}

\begin{equation}
\varepsilon=\frac{max\left|\lambda\left(\Lambda_{s}\right)\right|}{min\left|\lambda\left(\Lambda_{f}\right)\right|}\ll1.\label{eq:15_SMALL}
\end{equation}

In principle, the idea of the linear transformation is not new, see
e.g. \cite{OM1998}, but the principal point of the developed algorithm
concerns evaluation of this transformation. We have developed the
efficient and robust method that produces the best possible (to the
leading order) decomposition with respect to existing multiple-scales
hierarchy (see the attachment and \cite{BG2013,BGM2008,BM2009} for
more details).


\section{Singularly perturbed profiles}

In this section we extend the previous study taking in account transport
terms (diffusion and advection in general case). According to the
REDIM method the manifold is constructed by implementation of the
invariance condition \cite{F1988,RF1990,GKZ2004a,GKZ2004b} to the
system vector field extended by the transport terms \cite{BM2007,BM2009}.
As it was shown in e.g. \cite{BM2007} the advection/convection term
does not affect the state space of a reacting system, thus, one needs
to account only for the diffusion term. Hence, when an appropriate (with
respect to the decomposition of motions defined by the chemical source
term) coordinate system (u,v) is defined the simplest mathematical
model to the leading order is the following

\begin{equation}
\left\{ \begin{array}{c}
\frac{du(x,t)}{dt}=F_{s}\left(u(x,t),\, v(x,t)\right)+L_{1,x}(u(x,t),\, v(x,t))\\
\frac{dv(x,t)}{dt}=\frac{1}{\epsilon}F_{f}\left(u(x,t),\, v(x,t)\right)+L_{2,x}(u(x,t),\, v(x,t))
\end{array}\right.\label{eq:SPP1}
\end{equation}

Additionally, one can see from (\ref{eq:SPP1}) that our main assumption
treats the transport term as slow comparatively with the fast component
of the vector field.

The slow system evolution is then controlled by 
\[
u(x,t)=\left(u_{1}(x,t),...,\, u_{m_{s}}(x,t)\right),
\]

which are assumed to change slowly comparatively to the fast variables
\[
v(x,t)=\left(v_{1}(x,t),...,\, v_{m_{f}}(x,t)\right),\quad m_{s}+m_{f}=n.
\]

The diffusion terms are represented first by very general
and smooth differential operators $L_{1,x}(u(x,t),\, v(x,t))$, $L_{2,x}(u(x,t),\, v(x,t))$.
 Suppose that $u(x,t)$, $v(x,t)$ are smooth functions.

Initial data for the system \ref{eq:SPP1} are 
\begin{equation}
u(x,0)=u_{0}(x),v(x,0)=v_{0}(x).\label{eq:indata}
\end{equation}

Recall that functions $F_{s},F_{f}$ are of the same order. Then $\left|\left|\frac{dU}{dt}\right|\right|\sim O\left(1\right)$
while $\left|\left|\frac{dV}{dt}\right|\right|\sim O\left(\frac{1}{\varepsilon}\right)$.
We suppose also that operators $L_{1,x}(u(x,t),\, v(x,t)),L_{2,x}(u(x,t),\, v(x,t))$
has the same order as $F_{s},F_{f}$.

The zero approximation $S$ of the slow invariant manifold in the
phase space $(u,v)$ (the space of species) is represented in the
implicit form 
\[
F_{f}(u,v)=0
\]

The initial profile is $\varGamma_{0}(x):=(u_{0}(x),v_{0}(x));\, u_{0}(x)=u(x,0),\, v_{0}(x)=v(x,0)$.
Denote $\varGamma(x,t)$ a profile that is the solution of (\ref{eq:SPP1})
at time $t$ with the initial profile (initial data) $\varGamma_{0}(x)$
.

Now, we are going to use the singular perturbed structure of the system
(\ref{eq:SPP1}). We will follow formally the classical scheme for
singular perturbed systems of ODE. All profiles are surfaces in the
phase space $\mathrm{R^{n}}$. If initial data $\varGamma_{0}(x)$
do not belongs to the slow manifold it will be projected to slow manifold
along the fast subspaces parallel to the coordinate fast subspace
$\mathrm{R^{m_{f}}}$. As result we obtain the slow initial profile
$\varGamma_{0,s}(x)$ on the slow manifold $F_{f}\left(u,\, v\right)=0$.
For simplicity we can suppose that $\varGamma_{0}(x)$ itself belongs
to $S$. In the zero approximation the evolution of the slow initial
profile is due to the slow equation 
\[
\frac{du(x,t)}{dt}=F_{s}\left(u(x,t),\, v(x,t)\right)+L_{1,x}(u(x,t),\, v(x,t))
\]
under additional condition 
\begin{equation}
F_{f}\left(u(x,t),\, v(x,t)\right)=0,\label{eq:RE-0}
\end{equation}
i.e. the zero approximation $\varGamma_{0}(x,t)$ of $\varGamma(x,t)$
belongs to $S$ for all $t$.

The set $RM:=\cup_{t\in(0,\infty)}\varGamma(x,t)$ is called the reaction-diffusion
manifold (REDIM) and $RM_{0}:=\cup_{t\in(0,\infty)}\varGamma_{0}(x,t)$
is its zero approximation (for $\epsilon=0$).

If the dimension of the profile is equal to $s$ ($\dim\varGamma(x,t)=s$)
then $\dim RM=\dim RM_{0}=s+1$.

To obtain the first order approximation of the profile in the implicit
form , it is necessary to differentiate the equation $F_{f}(x(t),y(t),\epsilon)=0$
as a function of $t$ and by virtue of the system (\ref{eq:SPP1}),
the result is 
\[
0=\frac{\partial F_{f}}{\partial u}\frac{\partial u}{\partial t}+\frac{\partial F_{f}}{\partial v}\frac{\partial v}{\partial t}=\frac{\partial F_{f}}{\partial u}(F_{s}(u(x,t),v(x,t))+L_{1,x}(u(x,t),\, v(x,t))+
\]
\[
\frac{\partial F_{f}}{\partial v}(\frac{1}{\epsilon}F_{f}(u(x,t),v(x,t))+L_{2,x}(u(x,t),\, v(x,t))).
\]
We used the short notation $\frac{\partial F_{f}}{\partial u}$ for
$\frac{\partial F_{f}}{\partial u}(u(x,t),v(x,t))$ and, correspondingly,
for $\frac{\partial F_{f}}{\partial v}$. If $\det\frac{\partial F_{f}}{\partial v}\neq0$,
then this equation can be written in the more convenient form 
\[
\begin{array}{c}
0=F_{f}(u(x,t),v(x,t))+\varepsilon\left(\frac{\partial F_{f}}{\partial v}\right)^{-1}\frac{\partial F_{f}}{\partial u}F_{s}(u(x,t),v(x,t))+.\\
\epsilon\left(\frac{\partial F_{f}}{\partial v}\right)^{-1}\frac{\partial F_{f}}{\partial u}L_{1,x}(u(x,t),\, v(x,t))+\epsilon L_{2,x}(u(x,t),\, v(x,t)).
\end{array}
\]

We recover the zero order approximation, on setting $\epsilon=0$.
Denote the first approximation coefficient 
\[
h_{1}(u(x,t),v(x,t)):=\left(\frac{\partial F_{f}}{\partial v}\right)^{-1}\frac{\partial F}{\partial u}F_{s}(u(x,t),v(x,t))+
\]

\[
\left(\frac{\partial F_{f}}{\partial v}\right)^{-1}\frac{\partial F}{\partial u}L_{1,x}(u(x,t),\, v(x,t)))+L_{2,x}(u(x,t),\, v(x,t))
\]

The first part is the first approximation for ODE, the second part
is the transport correction. Thus, the first order approximation is
\begin{align}
0=F_{f}(u(x,t),v(x,t))+\epsilon h_{1}(u(x,t),v(x,t)).\label{for 5-1}
\end{align}
The transport terms $\varepsilon L_{1,x}(u(x,t),\, v(x,t)),L_{2,x}(u(x,t),\, v(x,t))$
have an influence only on the first approximation.

To shorten formulas we shall drop variable $(x,t)$ when it is clear
from the context.

To obtain the second order approximation it is necessary to differentiate
the equation $F_{f}(u(x,t),v(x,t))=0$ twice with respect to $t$
and equate it to zero:

\begin{equation}
\begin{array}{c}
0=\frac{d^{2}}{dt^{2}}(F_{f}(u,v))=\frac{d}{dt}\left(\frac{\partial F_{f}}{\partial u}\frac{\partial u}{\partial t}+\frac{\partial F_{f}}{\partial v}\frac{\partial v}{\partial t}\right)=\\
\left(\frac{\partial^{2}F_{f}}{\partial u^{2}}\frac{\partial u}{\partial t}+\frac{\partial^{2}F_{f}}{\partial u\partial v}\frac{\partial v}{\partial t}\right)\frac{\partial u}{\partial t}+\frac{\partial F_{f}}{\partial u}\frac{\partial^{2}u}{\partial t^{2}}\\
\left(\frac{\partial^{2}F_{f}}{\partial u\partial v}\frac{\partial u}{\partial t}+\frac{\partial^{2}F_{f}}{\partial v^{2}}\frac{\partial v}{\partial t}\right)\frac{\partial v}{\partial t}+\frac{\partial F_{f}}{\partial v}\frac{\partial^{2}v}{\partial t^{2}}
\end{array}\label{for 6-1-1}
\end{equation}

By using the system (\ref{eq:SPP1}) and linearity of operators $L_{1,x}(u(x,t),\, v(x,t)),\, L_{2,x}(u(x,t),\, v(x,t))$
we obtain 
\begin{equation}
\begin{array}{c}
\frac{\partial F_{f}}{\partial v}\frac{\partial^{2}v}{\partial t^{2}}=\frac{\partial F_{f}}{\partial v}\frac{\partial}{\partial t}\Big(\frac{1}{\epsilon}F_{f}(u,v)+L_{2,x}(u(x,t),\, v(x,t))\Big)=\\
\frac{\partial F_{f}}{\partial v}\Big(\frac{1}{\epsilon}\left[\frac{\partial F_{f}}{\partial u}\frac{\partial u}{\partial t}+\frac{\partial F_{f}}{\partial v}\frac{\partial v}{\partial t}\right]+L_{2,x}(\frac{\partial u(x,t)}{\partial t},\,\frac{\partial v(x,t)}{\partial t})\Big)=\\
\frac{\partial F_{f}}{\partial v}\Big(\frac{1}{\epsilon}\frac{\partial F_{f}}{\partial v}\left[\frac{1}{\epsilon}F_{f}(u,v)+L_{2,x}(u(x,t),\, v(x,t))\right]+\\
\frac{1}{\epsilon}\frac{\partial F_{f}}{\partial u}(F_{s}(u,v)+L_{1,x}(u(x,t),\, v(x,t))+L_{2,x}(\frac{\partial u(x,t)}{\partial t},\,\frac{\partial v(x,t)}{\partial t})\Big)\\
\frac{\partial F_{f}}{\partial u}\frac{\partial^{2}u}{\partial t^{2}}=\frac{\partial F_{f}}{\partial u}\frac{\partial}{\partial t}\Big(F_{s}(u,v)+L_{1,x}(u(x,t),\, v(x,t))\Big)=\\
\frac{\partial F_{f}}{\partial u}\Big(\frac{\partial F_{s}}{\partial u}\frac{\partial u}{\partial t}+\frac{\partial F_{s}}{\partial v}\frac{\partial v}{\partial t}+L_{1,x}(\frac{\partial u(x,t)}{\partial t},\,\frac{\partial v(x,t)}{\partial t})\Big)
\end{array}\label{for 7-1}
\end{equation}

Equation $(\ref{for 5})$ produces $\frac{1}{\epsilon}F(u,v)=-h_{1}(u,v)$,
and $\frac{\partial v}{\partial t}=\frac{1}{\epsilon}F_{f}(u,v)+L_{2,x}(u(x,t),\, v(x,t))$,
thus, equation $(\ref{for 7-1})$ can be rewritten as

\begin{equation}
\begin{array}{c}
\frac{\partial F_{f}}{\partial u}\frac{\partial^{2}u}{\partial t^{2}}=\frac{\partial F_{f}}{\partial u}\Big(\frac{\partial F_{s}}{\partial v}(L_{2,x}(u(x,t),\, v(x,t))-h_{1}(u,v))+\\
\frac{\partial F_{s}}{\partial u}(F_{s}(u,v)+L_{1,x}(u(x,t),\, v(x,t))+L_{1,x}(\frac{\partial u(x,t)}{\partial t},\,\frac{\partial v(x,t)}{\partial t})\Big)
\end{array}\label{for 9-1}
\end{equation}
Similarly, by using equation $(\ref{for 5})$ and equations $(\ref{for 1})$
\begin{equation}
\begin{array}{c}
\left(\frac{\partial^{2}F_{f}}{\partial u^{2}}\frac{\partial u}{\partial t}+\frac{\partial^{2}F_{f}}{\partial u\partial v}\right)\frac{\partial u}{\partial t}=\frac{\partial^{2}F_{f}}{\partial u^{2}}(\frac{\partial u}{\partial t})^{2}+\frac{\partial^{2}F_{f}}{\partial u\partial v}\frac{\partial v}{\partial t}\frac{\partial u}{\partial t}=\\
\frac{\partial^{2}F_{f}}{\partial u^{2}}(L_{1,x}(u(x,t),\, v(x,t))+F_{s}(u,v))^{2}+\frac{\partial^{2}F_{f}}{\partial u\partial v}(F_{s}(u,v)+\\
L_{1,x}(u(x,t),\, v(x,t)))(L_{2,x}(u(x,t),\, v(x,t))-h_{1}(u,v)),
\end{array}\label{for 10-1}
\end{equation}

\begin{equation}
\begin{array}{c}
\left(\frac{\partial^{2}F_{f}}{\partial u\partial v}\frac{\partial u}{\partial t}+\frac{\partial^{2}F_{f}}{\partial v^{2}}\frac{\partial v}{\partial t}\right)\frac{\partial v}{\partial t}=\frac{\partial^{2}F_{f}}{\partial v\partial u}\frac{\partial u}{\partial t}\frac{\partial v}{\partial t}+\frac{\partial^{2}F_{f}}{\partial v^{2}}(\frac{\partial v}{\partial t})^{2}=\\
\frac{\partial^{2}F_{f}}{\partial v\partial u}(L_{2,x}(u(x,t),\, v(x,t))-h_{1}(u,v))(F_{s}(u,v)+L_{1,x}(u(x,t),\, v(x,t)))+\\
\frac{\partial^{2}F}{\partial v^{2}}(L_{1,x}(u(x,t),\, v(x,t))-h_{1}(u,v))^{2}
\end{array}\label{for 11-1}
\end{equation}
Substitution of expressions $(\ref{for 7-1}),(\ref{for 9-1}),(\ref{for 10-1}),(\ref{for 11-1})$
in $(\ref{for 6-1-1})$ provides

\begin{equation}
\begin{array}{c}
0=\Big(\frac{1}{\epsilon}\Big)^{2}\Big(\frac{\partial F_{f}}{\partial v}\Big)^{2}F_{f}(u,v)+\frac{1}{\epsilon}\Big(\frac{\partial F_{f}}{\partial v}\Big)^{2}L_{2,x}(u(x,t)\, v(x,t))+\frac{1}{\epsilon}\frac{\partial F_{f}}{\partial u}\frac{\partial F_{f}}{\partial v}(F_{s}(u,v)+\\
L_{1,x}(u(x,t),\, v(x,t)))+\frac{\partial F_{f}}{\partial v}L_{1,x}(\frac{\partial u(x,t)}{\partial t},\,\frac{\partial v(x,t)}{\partial t}))+\\
\frac{\partial F_{f}}{\partial u}\Big(\frac{\partial F_{s}}{\partial v}(L_{2,x}(u(x,t),\, v(x,t)))-h_{1}(u,v))+\frac{\partial F_{s}}{\partial u}(F_{s}(u,v)+\\
L_{1,x}(u(x,t),\, v(x,t)))+L_{2,x}(\frac{\partial u(x,t)}{\partial t},\,\frac{\partial v(x,t)}{\partial t})\Big)+\\
\frac{\partial^{2}F_{s}}{\partial v^{2}}(L_{2,x}(u(x,t),\, v(x,t)))-h_{1}(u,v))^{2}+\frac{\partial^{2}F_{f}}{\partial u\partial v}(F_{s}(u,v)+\\
L_{1,x}(u(x,t),\, v(x,t)))(L_{2,x}(u(x,t),\, v(x,t))-h_{1}(u,v))+\\
\frac{\partial^{2}F_{f}}{\partial v\partial u}(L_{2,x}(u(x,t),\, v(x,t)))-h_{1}(u,v))(F_{s}(u,v)+\\
L_{1,x}(u(x,t),\, v(x,t))))+\frac{\partial^{2}F_{f}}{\partial u^{2}}(L_{1,x}(u(x,t),\, v(x,t))+F_{s}(u,v))^{2}.
\end{array}\label{for 12-1}
\end{equation}
Multiplication of $(\ref{for 12-1})$ by $\epsilon^{2}(\frac{\partial F_{f}}{\partial v})^{-2}$
yields 
\begin{equation}
\begin{array}{c}
0=F_{f}(u,v)+\epsilon L_{2,x}(u(x,t),\, v(x,t))+\epsilon\Big(\frac{\partial F_{f}}{\partial v}\Big)^{-1}\frac{\partial F_{f}}{\partial u}(F_{s}(u,v)+L_{1,x}(u(x,t),\, v(x,t)))+\\
\epsilon^{2}\Big(\frac{\partial F_{f}}{\partial v}\Big)^{-2}\Big(\frac{\partial F_{f}}{\partial v}L_{2,x}(\frac{\partial u(x,t)}{\partial t},\,\frac{\partial v(x,t))}{\partial t})+\frac{\partial F_{f}}{\partial u}\Big(\frac{\partial F_{s}}{\partial v}(L_{2,x}(u(x,t),\, v(x,t)))-h_{1}(u,v))+\\
\frac{\partial F_{s}}{\partial u}((F_{s}(u,v)+L_{1,x}(u(x,t),\, v(x,t)))+L_{1,x}(\frac{\partial u(x,t)}{},\,\frac{\partial v}{\partial t}))\Big)+\\
\frac{\partial^{2}F_{f}}{\partial v^{2}}(L_{2,x}(u(x,t),\, v(x,t)))-h_{1}(u,v))^{2}+2\frac{\partial^{2}F_{f}}{\partial u\partial v}(F_{s}(u,v)+\\
L_{1,x}(u(x,t),\, v(x,t)))(L_{2,x}(u(x,t),\, v(x,t))-h_{1}(u,v))+\\
\frac{\partial^{2}F_{f}}{\partial u^{2}}(F_{s}(u,v)+L_{1,x}(u(x,t),\, v(x,t)))^{2}\Big).
\end{array}\label{for 13-1}
\end{equation}
Denote 
\begin{equation}
\begin{array}{c}
h_{2}(u,v):=\Big(\frac{\partial F_{f}}{\partial v}\Big)^{-2}\Big(\frac{\partial F_{f}}{\partial v}L_{2,x}(\frac{\partial u(x,t)}{\partial t},\,\frac{\partial v(x,t))}{\partial t})+\frac{\partial F_{f}}{\partial u}\Big(\frac{\partial F_{s}}{\partial v}(L_{2,x}(u(x,t),\, v(x,t)))-h_{1}(u,v))+\\
\frac{\partial F_{s}}{\partial u}((F_{s}(u,v)+L_{1,x}(u(x,t),\, v(x,t))))+L_{1,x}(\frac{\partial u(x,t)}{\partial t},\,\frac{\partial v(x,t)}{\partial t}))\Big)+\\
\frac{\partial^{2}F_{f}}{\partial v^{2}}(L_{2,x}(u(x,t),\, v(x,t)))-h_{1}(u,v))^{2}+2\frac{\partial^{2}F_{f}}{\partial u\partial v}(F_{s}(u,v)+\\
L_{1,x}(u(x,t),\, v(x,t)))(L_{2,x}(u(x,t),\, v(x,t))-h_{1}(u,v))+\\
\frac{\partial^{2}F_{f}}{\partial u^{2}}(F_{s}(u,v)+L_{1,x}(u(x,t),\, v(x,t)))^{2}\Big),\\
\\
\end{array}\label{for 14-1}
\end{equation}
therefore, second order approximation is given by 
\begin{align}
0=F(u,v)+\epsilon h_{1}(u,v)+\epsilon^{2}h_{2}(u,v).\label{for 15-1}
\end{align}

\subsection{Slow transport term}

The special case that is typical for chemical kinetics models is the
case of relatively slow transport processes, namely,

\begin{equation}
\left\{ \begin{array}{c}
\frac{du(x,t)}{dt}=F_{s}\left(u(x,t),\, v(x,t)\right)+\varepsilon L_{1,x}(u(x,t),\, v(x,t))\\
\varepsilon\frac{dv(x,t)}{dt}=F_{f}\left(u(x,t),\, v(x,t)\right)+\varepsilon L_{2,x}(u(x,t),\, v(x,t))
\end{array}\right.\label{eq:SPP2}
\end{equation}

For this special case evolution of the slow initial profile for $\varepsilon=0$
is due to the slow equation 
\[
\frac{du(x,t)}{dt}=F_{s}\left(u(x,t),\, v(x,t)\right),
\]
under additional condition 
\[
F_{f}\left(u(x,t),\, v(x,t)\right)=0,
\]
i.e. on the slow manifold $F_{f}\left(u,\, v\right)=0.$

The transport term $\varepsilon L_{1,x}(U(x,t),\, V(x,t))$ has an
influence only on the first order approximation. In this special case
evolution of any point of the initial profile coincides with the evolution
of corresponding trajectory on the slow manifold.

\subsection{First order approximation of the slow invariant manifold for Laplacian}

Consider the following Singularly Perturbed System of equations, where
the transport operator is given by the Laplace operator: 
\begin{align}
 & \frac{\partial u}{\partial t}=\frac{1}{\epsilon}F(u,v)+\Delta u\label{for 1}\\
 & \frac{\partial v}{\partial t}=G(u,v)+\Delta v
\end{align}
where $0<\epsilon<<1$, $u=u(x,t),v=v(x,t)$. System $(\ref{for 1})$
can be reformulated as 
\begin{align}
 & \epsilon\frac{\partial u}{\partial t}=F(u,v)+\epsilon\Delta u\label{for2}\\
 & \frac{\partial v}{\partial t}=G(u,v)+\Delta v
\end{align}

We will find an implicit form of the first approximation of the slow
invariant manifold slow manifold

\begin{align*}
 & F(u,v)=0
\end{align*}
First approximation is obtained by differentiation of $F(u,v)=0$
with respect to $t$ 
\begin{align}
0=\frac{\partial F}{\partial u}\frac{\partial u}{\partial t}+\frac{\partial F}{\partial v}\frac{\partial v}{\partial t}=\frac{\partial F}{\partial u}(\frac{1}{\epsilon}F(u,v)+\Delta u)+\frac{\partial F}{\partial v}(G(u,v)+\Delta v)\label{for 3}
\end{align}
Multiplication of $(\ref{for 3})$ by $\epsilon(\frac{\partial F}{\partial u})^{-1}$
yields, 
\begin{align}
0=F(u,v)+\varepsilon\left[(\frac{\partial F}{\partial u})^{-1}\frac{\partial F}{\partial v}G(u,v)+(\frac{\partial F}{\partial u})^{-1}\frac{\partial F}{\partial v}G(u,v)\Delta v+\Delta u\right]\label{for 4}
\end{align}
Denote the first approximation coefficient 
\begin{align}
h_{1}(u,v)=(\frac{\partial F}{\partial u})^{-1}\frac{\partial F}{\partial v}G(u,v)+\left[(\frac{\partial F}{\partial u})^{-1}\frac{\partial F}{\partial v}G(u,v)\Delta v+\Delta u\right]
\end{align}
The first part is the first approximation for ODE, the second part
is the transport correction.

Thus, the first order approximation is given by 
\begin{align}
0=F(u,v)+\epsilon h_{1}(u,v)\label{for 5}
\end{align}

The first approximation for slow transport terms \ref{eq:SPP2} is

\[
h_{1}(u,v)=(\frac{\partial F}{\partial u})^{-1}\frac{\partial F}{\partial v}G(u,v)+\left[(\frac{\partial F}{\partial u})^{-1}\frac{\partial F}{\partial v}G(u,v)\Delta v\right]
\]

\subsection{Second order approximation of the slow manifold for Laplacian}

In order to obtain the second order approximation it is necessary
to differentiate the equation $F(u,v)=0$ twice with respect to $t$
\begin{align}
 & 0=\frac{d^{2}}{dt^{2}}(F(u,v))=\frac{d}{dt}(\frac{\partial F}{\partial u}\frac{\partial u}{\partial t}+\frac{\partial F}{\partial v}\frac{\partial v}{\partial t})=(\frac{\partial^{2}F}{\partial u^{2}}\frac{\partial u}{\partial t}+\frac{\partial^{2}F}{\partial u\partial v}\frac{\partial v}{\partial t})\frac{\partial u}{\partial t}+\\
 & \frac{\partial F}{\partial u}\frac{\partial^{2}u}{\partial t^{2}}+(\frac{\partial^{2}F}{\partial v\partial u}\frac{\partial u}{\partial t}+\frac{\partial^{2}F}{\partial v^{2}}\frac{\partial v}{\partial t})\frac{\partial v}{\partial t}+\frac{\partial F}{\partial v}\frac{\partial^{2}v}{\partial t^{2}}\label{for 6}
\end{align}

where 
\begin{align}
 & \frac{\partial F}{\partial u}\frac{\partial^{2}u}{\partial t^{2}}=\frac{\partial F}{\partial u}\frac{\partial}{\partial t}\Big(\frac{1}{\epsilon}F(u,v)+\Delta u\Big)=\frac{\partial F}{\partial u}\Big(\frac{1}{\epsilon}(\frac{\partial F}{\partial u}\frac{\partial u}{\partial t}+\frac{\partial F}{\partial v}\frac{\partial v}{\partial t})+\Delta(\frac{\partial u}{\partial t})\Big)=\\
 & \frac{\partial F}{\partial u}\Big(\frac{1}{\epsilon}\frac{\partial F}{\partial u}(\frac{1}{\epsilon}F(u,v)+\Delta u)+\frac{1}{\epsilon}\frac{\partial F}{\partial v}(G(u,v)+\Delta v)+\Delta(\frac{\partial u}{\partial t})\Big)\label{for 7}
\end{align}
\begin{align}
\frac{\partial F}{\partial v}\frac{\partial^{2}v}{\partial t^{2}}=\frac{\partial F}{\partial v}\frac{\partial}{\partial t}\Big(G(u,v)+\Delta v\Big)=\frac{\partial F}{\partial v}\Big(\frac{\partial G}{\partial u}\frac{\partial u}{\partial t}+\frac{\partial G}{\partial v}\frac{\partial v}{\partial t}+\Delta(\frac{\partial v}{\partial t})\Big).\label{for 8}
\end{align}
Equation $(\ref{for 5})$ produces, $\frac{1}{\epsilon}F(u,v)=-h_{1}(u,v)$,
and $\frac{\partial u}{\partial t}=\frac{1}{\epsilon}F(u,v)+\Delta u$,
thus, the equation $(\ref{for 8})$ can be rewritten as 
\begin{align}
\frac{\partial F}{\partial v}\frac{\partial^{2}v}{\partial t^{2}}=\frac{\partial F}{\partial v}\Big(\frac{\partial G}{\partial u}(\Delta u-h_{1}(u,v))+\frac{\partial G}{\partial v}(G(u,v)+\Delta v)+\Delta(\frac{\partial v}{\partial t})\Big)\label{for 9}
\end{align}
Similarly, by using equation $(\ref{for 5})$ and equations $(\ref{for 1})$
\begin{align}
 & (\frac{\partial^{2}F}{\partial u^{2}}\frac{\partial u}{\partial t}+\frac{\partial^{2}F}{\partial u\partial v}\frac{\partial v}{\partial t})\frac{\partial u}{\partial t}=\frac{\partial^{2}F}{\partial u^{2}}(\frac{\partial u}{\partial t})^{2}+\frac{\partial^{2}F}{\partial u\partial v}\frac{\partial v}{\partial t}\frac{\partial u}{\partial t}=\\
 & \frac{\partial^{2}F}{\partial u^{2}}(\Delta u-h_{1}(u,v))^{2}+\frac{\partial^{2}F}{\partial u\partial v}(G(u,v)+\Delta v)(\Delta u-h_{1}(u,v))\label{for 10}
\end{align}
\begin{align}
 & (\frac{\partial^{2}F}{\partial v\partial u}\frac{\partial u}{\partial t}+\frac{\partial^{2}F}{\partial v^{2}}\frac{\partial v}{\partial t})\frac{\partial v}{\partial t}=\frac{\partial^{2}F}{\partial v\partial u}\frac{\partial u}{\partial t}\frac{\partial v}{\partial t}+\frac{\partial^{2}F}{\partial v^{2}}(\frac{\partial v}{\partial t})^{2}=\\
 & \frac{\partial^{2}F}{\partial v\partial u}(\Delta u-h_{1}(u,v))(G(u,v)+\Delta v)+\frac{\partial^{2}F}{\partial v^{2}}(G(u,v)+\Delta v)^{2}\label{for 11}
\end{align}
Substitution of expressions $(\ref{for 8}),(\ref{for 9}),(\ref{for 10}),(\ref{for 11})$
in $(\ref{for 6})$ provides 
\begin{align}
 & 0=\Big(\frac{1}{\epsilon}\Big)^{2}\Big(\frac{\partial F}{\partial u}\Big)^{2}F(u,v)+\frac{1}{\epsilon}\Big(\frac{\partial F}{\partial u}\Big)^{2}\Delta u+\frac{1}{\epsilon}\frac{\partial F}{\partial u}\frac{\partial F}{\partial v}(G(u,v)+\Delta v)+\frac{\partial F}{\partial u}\Delta(\frac{\partial u}{\partial t})+\\
 & \frac{\partial F}{\partial v}\Big(\frac{\partial G}{\partial u}(\Delta u-h_{1}(u,v))+\frac{\partial G}{\partial v}(G(u,v)+\Delta v)+\Delta(\frac{\partial v}{\partial t})\Big)+\\
 & \frac{\partial^{2}F}{\partial u^{2}}(\Delta u-h_{1}(u,v))^{2}+\frac{\partial^{2}F}{{\partial u}{\partial v}}(G(u,v)+\Delta v)(\Delta u-h_{1}(u,v))+\\
 & \frac{\partial^{2}F}{\partial v\partial u}(\Delta u-h_{1}(u,v))(G(u,v)+\Delta v)+\frac{\partial^{2}F}{\partial v^{2}}(G(u,v)+\Delta v)^{2}\label{for 12}
\end{align}
Multiplication of $(\ref{for 12})$ by $\epsilon^{2}(\frac{\partial F}{\partial u})^{-2}$
yields 
\begin{align}
 & 0=F(u,v)+\epsilon\Delta u+\epsilon\Big(\frac{\partial F}{\partial u}\Big)^{-1}\frac{\partial F}{\partial v}(G(u,v)+\Delta v)+\\
 & \epsilon^{2}\Big(\frac{\partial F}{\partial u}\Big)^{-2}\Big(\frac{\partial F}{\partial u}\Delta(\frac{\partial u}{\partial t})+\frac{\partial F}{\partial v}\Big(\frac{\partial G}{\partial u}(\Delta u-h_{1}(u,v))+\frac{\partial G}{\partial v}(G(u,v)+\Delta v)+\Delta(\frac{\partial v}{\partial t})\Big)+\\
 & \frac{\partial^{2}F}{\partial u^{2}}(\Delta u-h_{1}(u,v))^{2}+2\frac{\partial^{2}F}{\partial u\partial v}(G(u,v)+\Delta v)(\Delta u-h_{1}(u,v))+\frac{\partial^{2}F}{\partial v^{2}}(G(u,v)+\Delta v)^{2}\Big).\label{for 13}
\end{align}
Finally, denote 
\begin{align}
 & h_{2}(u,v)=\Big(\frac{\partial F}{\partial u}\Big)^{-2}\Big(\frac{\partial F}{\partial u}\Delta(\frac{\partial u}{\partial t})+\frac{\partial F}{\partial v}\Big(\frac{\partial G}{\partial u}(\Delta u-h_{1}(u,v))+\frac{\partial G}{\partial v}(G(u,v)+\Delta v)+\Delta(\frac{\partial v}{\partial t})\Big)+\\
 & \frac{\partial^{2}F}{\partial u^{2}}(\Delta u-h_{1}(u,v))^{2}+2\frac{\partial^{2}F}{\partial u\partial v}(G(u,v)+\Delta v)(\Delta u-h_{1}(u,v))+\frac{\partial^{2}F}{\partial v^{2}}(G(u,v)+\Delta v)^{2}\Big)\label{for 14}
\end{align}
Therefore, the second order approximation is given by 
\begin{align}
0=F(u,v)+\epsilon h_{1}(u,v)+\epsilon^{2}h_{2}(u,v).\label{for 15}
\end{align}

\begin{figure}[h]
 \centering
 \includegraphics{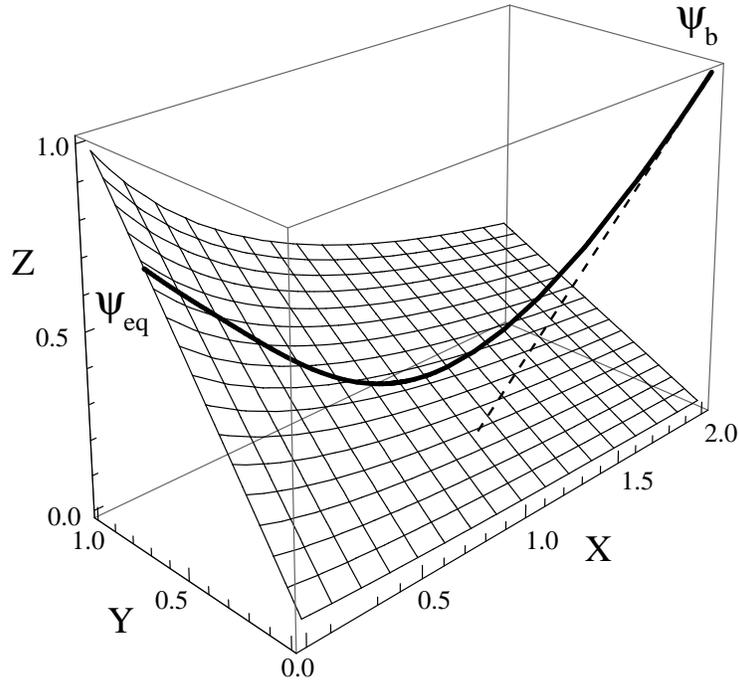}
 \caption{System state space (X,Y,Z) with 2D slow manifold (mesh) and stationary solution profile (black solid curve).
 The approximation of the fast part of the system solution trajectory is shown by dashed line.}
 \label{fig:F1}
 \end{figure}

\section{3D Michaelis-Menten model with Laplacian}

An illustration of the suggested approximation scheme for a system
with the diffusion is presented for a benchmark model below. We consider
the 3D Michaelis-Menten model. The mathematical model consists of
three ODEs

\begin{align}
 & \frac{dX}{dt}=-XZ+L_{1}(1-Z-\mu(1-Y))\\
 & \frac{dY}{dt}=-L_{3}YZ+\frac{L_{4}}{L_{2}}(1-Y)\\
 & \frac{dZ}{dt}=\frac{1}{L_{2}}((-XZ+1-Z-\mu(1-Y))+\mu)-L_{3}YZ+\frac{L_{4}}{L_{2}}(1-Y)))\label{for 16}
\end{align}

The system parameters are $L_{1}=0.99$, $L_{2}=1$, $L_{3}=0.05$,
$L_{4}=0.1$, $\mu=1$ (see e.g. \cite{RF2001} for details and references).
By taking the diffusion into account we obtain the following PDE
system with the diffusion constant $\delta=0.01$:

\begin{align}
 & \frac{\partial X}{\partial t}=-XZ+L_{1}(1-Z-\mu(1-Y))+\delta\Delta X\label{eq: for 15}\\
 & \frac{\partial Y}{\partial t}=-L_{3}YZ+\frac{L_{4}}{L_{2}}(1-Y)+\delta\Delta Y\\
 & \frac{\partial Z}{\partial t}=\frac{1}{L_{2}}((-XZ+1-Z-\mu(1-Y))+\mu)-L_{3}YZ+\frac{L_{4}}{L_{2}}(1-Y)))+\delta\Delta Z\label{for 17}
\end{align}


The system (\ref{eq: for 15}) - (\ref{for 17}) is considered
with the following initial and boundary conditions:

\begin{align}
 & \begin{pmatrix} 
 X(t,0)=X_{eq} \\
 Y(t,0)=Y_{eq} \\
 Z(t,0)=Z_{eq}\\
 \end{pmatrix}\label{eq: for 18}\\
 & \begin{pmatrix} 
 X(t,1)=2 \\
 Y(t,1)=0 \\
 Z(t,1)=1\\
 \end{pmatrix}\\
 & \begin{pmatrix} 
 X(0,x)=(2-X_{eq})x+X_{eq} \\
 Y(0,x)=(-Y_{eq})x+Y_{eq} \\
 Z(0,x)=(1-Z_{eq})x+Z_{eq}\\
 \end{pmatrix}\label{for 20}
\end{align}

Here $(X_{eq},Y_{eq},Z_{eq})$ are coordinates of the equilibrium point that are calculated numerically and $x$ is spatial variable.
Initial conditions are chosen to be straight lines, they satisfy the general assumption
- join initial and equilibrium values on boundaries. 

As in the previous section the main assumption is the transport term
is slow compared with the fast vector field. First, several numerical
experiments were performed drawing a 2D slow manifold and stationary
solution of the system. Figure \ref{fig:F1} shows a connection between
the zero approximation of the slow manifold and the profile of the
stationary system solution to the PDE in the original coordinates
$(x,y,z)$. In the Fig. \ref{fig:F1} the stationary solution can
be roughly divided into two parts: the slow part of the stationary
solution that is very close to the slow invariant manifold and the
second one, which is asymptotically close to the fast sub-field
of the system shown by the dashed line. The small third part conjugates
the slow and fast parts and can be asymptotically neglected. This simulation
demonstrates that the proposed machinery of singularly perturbed profiles
is relevant to the problem under consideration. It means that the
original system (\ref{eq: for 15}) - (\ref{for 17}) can be decomposed
into two subsystems: the slow one and the fast one.

In order to implement the scheme described above we apply the GQL
method \cite{BGM2008,BGM2011} to convert the model to a system in
the standard SPS form in the new coordinate set (U,V,W). 
The transformation between $(X,Y,Z)$ and (U,V,W) is 

\begin{align}
 &X=0.73U-0.19V-0.66W\label{eq: for 21}\\
 & Y=-0.05U+0.94V-0.33W\\
 & Z=0.68U+0.28V+0.68W.\label{for 22}
\end{align}

 \begin{figure}[h]
 \centering
 \includegraphics{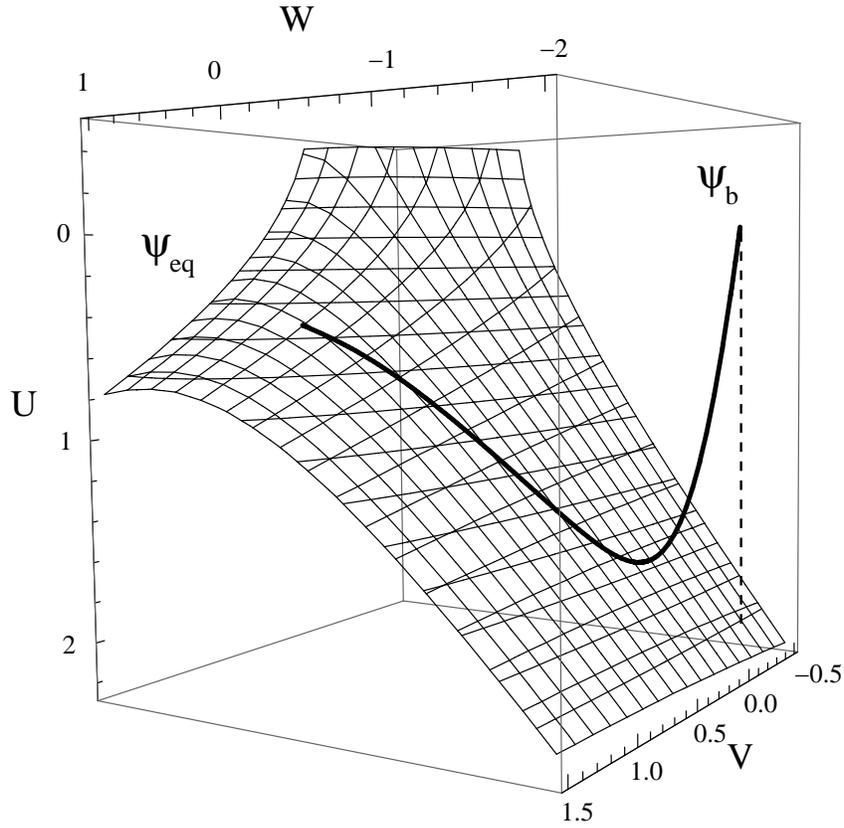}
 \caption{System state space (U,V,W) after the coordinate transformation (\ref{eq: for 21}-\ref{for 22}) into the standard SPS form, 2D slow manifold (mesh) and stationary solution profile (black solid curve). The approximation of the fast part of the system solution trajectory is shown by dashed line.}
 \label{fig:F2}
 \end{figure}

Figure \ref{fig:F2} shows the zero approximation of the slow manifold
and trajectory to PDE`s solution in the decomposed coordinate system
\cite{BG2008}. As it has to be in Fig. \ref{fig:F2} the fast part (which does
not belong to the manifold and shown by a dashed line) is approximated by
a parallel line to $U$ coordinate.

 \begin{figure}[h]
 \centering
 \includegraphics{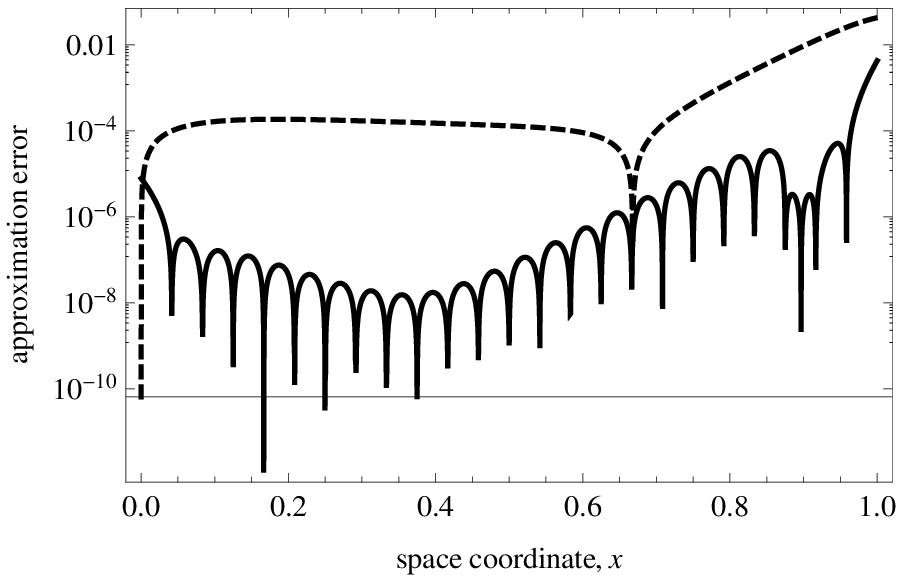}
 \caption{Dashed line matches  the solution profile on the zero order approximation and thick line matches the solution profile on the first order approximation with a logarithmically scaled y axis.}
 \label{fig:F3}
 \end{figure}

At present we concentrate mainly on the slow part, because the fast
part is very simple and asymptotically corresponds to one dimensional
fast manifold defining a projection of the initial data to the slow
manifold in pure homogeneous system. Thus, following the main analytic
results we compare the zero and the first order approximation
of the slow part of system solution stationary profiles.
Figure \ref{fig:F3} illustrates the stationary solution profiles on the
zero order approximation 
\[
H_0:=F(U(x,t),V(x,t),W(x,t))=0
\]
and on the first order approximation
\[ 
H_1:=F(U(x,t),V(x,t),W(x,t))+\epsilon h_1(U(x,t),V(x,t),W(x,t))=0
\]

of the system slow manifold.

Here $H_0:=F(U,V,W)$ is an implicit form of 2D slow manifold for 3D Michaelis-Menten model given by
the left-hand side of the fast part (for variable $U$) of the SPS system, namely,

\[
H_0=F(U,V,W)=0.007(0.06-0.7U^2+0.07V^2+U(-1.03-0.12V-0.06W)+\\
 \]
 \[V(0.88+0.42W)-1.4W+0.63W^2)=0
\] 

Figure \ref{fig:F3} shows $H_0(U(x,t),V(x,t),W(x,t))$ and $H_1(U(x,t),V(x,t),W(x,t))$ approximations given in the implicit form. Thus, the error of the approximation can be easily estimated. In this figure we can see that the first order approximation of the profile is more
accurate than the zero order approximation.

 \begin{figure}[h]
 \centering
 \includegraphics{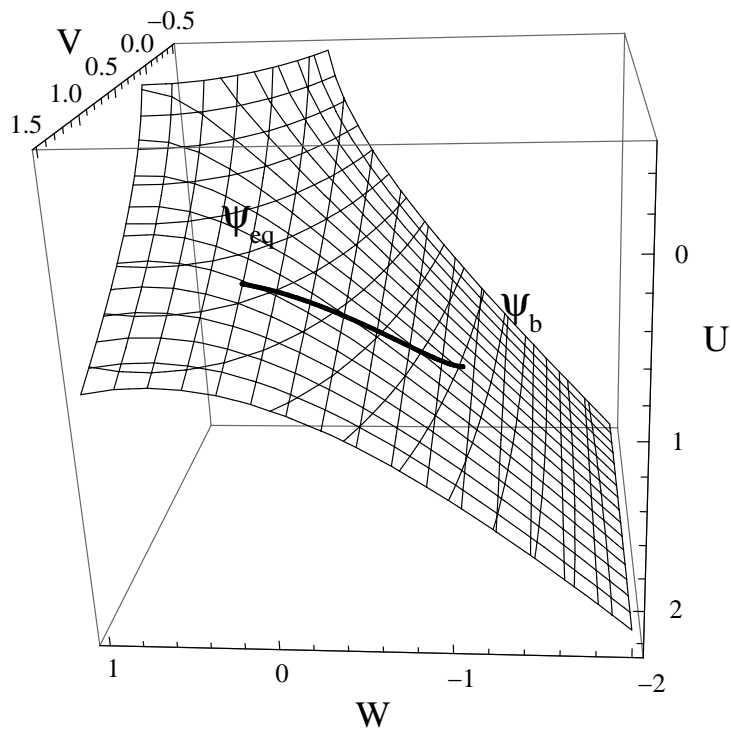}
 \caption{System state space (U,V,W) after transformation to the standard SPS
form, 2D slow manifold (mesh) and stationary solution profile (black curve)
for model with boundary condition which is very close to slow manifold.}
 \label{fig:F4}
 \end{figure}

The next Fig. \ref{fig:F4} illustrates  the zero order approximation of the slow manifold
and the solution profile trajectory of the PDE`s solution in the decomposed coordinate system
for the system with boundary condition close to the slow manifold. In Fig. \ref{fig:F4} the
fast part of the system solution profile, which does not belong to the manifold, vanishes.

\begin{figure}[h]
 \centering
 \includegraphics{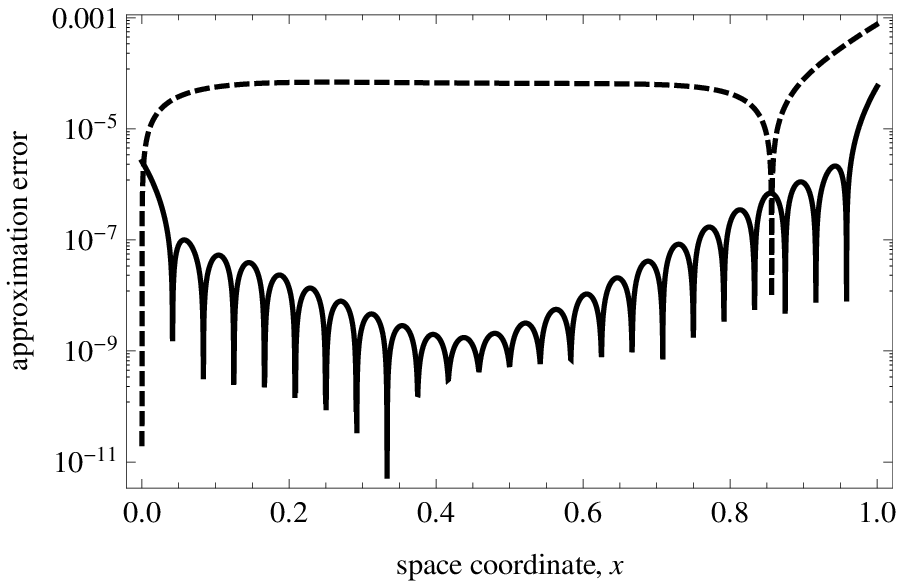}
 \caption{The system with boundary condition close to slow manifold. Dashed line matches  the solution profile on the zero order approximation and thick line matches the solution profile on the first order approximation with a logarithmically scaled y axis.}
 \label{fig:F5}
 \end{figure}
 
Similarly to Fig. \ref{fig:F3}, Fig. \ref{fig:F5} illustrates the stationary solution profiles
on the zero order approximation given in the implicit form by the RHS of U, namely, by
\[
H_0(U(x,t),V(x,t),W(x,t)):=F(U(x,t),V(x,t),W(x,t))
\]
and by the first order correction term correspondingly $H_1:=(U(x,t),V(x,t),W(x,t)):$
$$F(U(x,t),V(x,t),W(x,t))+\epsilon h_1(U(x,t),V(x,t),W(x,t))$$
By comparing of Figs. \ref{fig:F3} and \ref{fig:F5} we can see that the approximation error for the system with the boundary condition placed closer to slow manifold is
less than for the system with the boundary condition far from slow manifold.

\section{Conclusion}

The formal functional series for the system in the standard SP form
was used to study and approximate the REDIM (as a manifold made of
the system solutions profiles). If the transport/diffusion term is
considered as a first order correction to the system in the standard
SPS form, under assumption we do know the transformation leading to
this special form, then the correction terms were estimated until
the second order in general form of the transport term. The spacial
form of the diffusion described by the Laplacian was used to illustrate
and estimate the correction terms explicitly. It was found that the
correction terms are spatially dependent and can be evaluated on e.g.
system solution profiles (stationary or transient). The results were
illustrated and validated by the well known example of Michaelis-Menten
enzyme model extended by the Laplace operator. Effects of the system
decomposition, influence of the boundary conditions and diffusion
transport was studied.

\section*{Acknowledgments}

The financial support by DFG, especially within the GIF (Grant 1162-148.6/2011)
project is gratefully acknowledged.

\end{document}